\documentclass[12pt]{amsart}
\usepackage{amsmath,amsthm,amsfonts,amssymb,eucal}

\renewcommand {\a}{ \alpha }

\newcommand{\vare}{\varepsilon}
\newcommand{\g}{\gamma}
\newcommand{\G}{\Gamma}

\newcommand{\s}{\sigma}
\renewcommand{\l}{\lambda}

\newcommand{\R}{ \mathbb R}

\newcommand{\N}{ \mathbb N}

\newcommand {\GH}{\mathfrak H}
\newcommand {\GF}{\mathfrak F}
\newcommand {\ba}{\mathbf a}
\newcommand {\bb}{\mathbf b}
\newcommand {\BA}{\mathbf A}

\newcommand {\BH}{\mathbf H}
\newcommand {\BI}{\mathbf I}
\newcommand {\BJ}{\mathbf J}
\newcommand {\BT}{\mathbf T}
\newcommand {\bh}{\mathbf h}

\newcommand{\CB}{\mathcal B}
\newcommand{\CA}{\mathcal A}

\newcommand{\CF}{\mathcal F}

\newcommand{\tH}{\textsf H}
\newcommand{\tL}{\textsf L}

\newcommand{\wt}{\widetilde}

\DeclareMathOperator{\ess}{ess} \DeclareMathOperator{\ac}{ac}
 \DeclareMathOperator {\re}
{{Re}} 
\DeclareMathOperator{\dom}{Dom}

\newtheorem{thm}{Theorem}[section]

\newtheorem{lem}[thm]{Lemma}
\newtheorem{prop}[thm]{Proposition}

\theoremstyle{definition}

\theoremstyle{remark}

\numberwithin{equation}{section}
\newcommand{\thmref}[1]{Theorem~\ref{#1}}

\begin{document}

\title[On the discrete spectrum of a family of operators]
{On the discrete spectrum of a family of differential operators}
\author[M. Solomyak]{Michael Solomyak}
\address{Department of Mathematics\\The Weizmann Institute of Science\\
Rehovot 76100\\Israel}
\email{michail.solomyak@weizmann.ac.il}
\date{15.01.2004}
\dedicatory {Dedicated to Victor Borisovich Lidskii on the occasion
of his 80-th birthday}

\begin{abstract}
A family $\BA_\a$ of differential operators depending on a real
parameter $\a$ is considered. The problem can be formulated in the
language of perturbation theory of quadratic forms. The
perturbation is only relatively bounded but not relatively compact
with respect to the unperturbed form.

The spectral properties of the operator $\BA_\a$ strongly depend
on $\a$. In particular, for $\a<\sqrt2$ the spectrum of $\BA_\a$
below $1/2$ is finite, while for $\a>\sqrt2$ the operator has no
eigenvalues at all. We study the asymptotic behaviour of the
number of eigenvalues as $\a\nearrow\sqrt2$. We reduce this
problem to the one on the spectral asymptotics for a certain
Jacobi matrix.

\end{abstract}
\maketitle

\section{Introduction}
In this paper we study the discrete spectrum of a family $\BA_\a$ of
differential operators in the space $\tL^2(\R^2)$, defined by the
differential expression
\begin{equation}\label{eq}
\CA U=-U''_{xx}+\frac1{2}\bigl(-U''_{yy}+y^2 U\bigr)
\end{equation}
and the ``transmission condition'' on the line $x=0$:
\begin{equation}\label{tran}
U'_x(+0,y)-U'_x(-0,y)=\a y\,U(0,y),\qquad y\in\R.
\end{equation}
In \eqref{tran} $\a$ is a real parameter (the coupling constant).
So, the differential expression which defines the action of the
operator does not involve $\a$. The parameter appears only in the
condition \eqref{tran} which defines the operator domain of
$\BA_\a$. The replacement $\a\mapsto -\a$ corresponds to the
change of variables $y\mapsto -y$ which does not affect the
spectrum. For this reason, below we discuss only $\a>0$.

As we shall see, the spectrum $\s(\BA_\a)$ of $\BA_\a$ has the
discrete component only for $\a<\sqrt2$, and we study its
behaviour as $\a\nearrow\sqrt2$.

\vskip0.2cm Operators, containing the family $\BA_\a$ as a special
case, were suggested by U. Smilansky \cite{sm} as a model of an
irreversible quantum system. Some important conclusions on the
spectrum of $\BA_\a$ were made in \cite{sm} ``on the physical
level of rigour''. The first mathematical results on the subject
were obtained in \cite{sol}. In \cite{nabsol} they are
considerably extended. The work on \cite{nabsol} is yet unfinished
and its results are not used in the present paper. However, below
we mention some of these results in the course of our general
discussion.

It was shown in \cite{sol} for $\a\neq\sqrt 2$ and in \cite{nabsol} for
 $\a=\sqrt 2$, that for any $\a\in\R$ the operator $\BA_\a$, defined
originally on the set of all functions from the Schwartz class,
satisfying the condition \eqref{tran}, admits a unique
self-adjoint realization. The special role of the value $\a=\sqrt
2$ will be explained later. Below we refer to the values $\a<\sqrt
2$ as {\sl small} and to the values $\a>\sqrt 2$ as {\sl large}.
The spectral properties of the operator $\BA_\a$ for the small and
the large values of the parameter $\a$ are quite different.

\vskip0.2cm

It is useful to consider the quadratic form $\ba_\a[U]$ which formally
corresponds to the operator $\BA_\a$. It can be written as
\begin{equation}\label{qform}
\ba_\a[U]=\ba_0[U]+\a\bb[U]
\end{equation}
where
\begin{gather}
\ba_0[U]=\int_{\R^2}\bigl(|U'_x|^2+\frac{1}{2}(|U'_y|^2+y^2|U|^2)\bigr)dxdy;
\label{ba0}\\ \bb[U]=\int_\R y|U(0,y)|^2dy.\label{bb}
\end{gather}
We view $\ba_0[U]$ as the unperturbed quadratic form and $\a\bb[U]$ as the
perturbation.
An important feature of the problem studied stems from the fact that
$\bb[U]$ is only relatively bounded but not relatively compact with respect
to the quadratic form $\ba_0[U]$. For this reason, the standard results of the
perturbation theory do not apply, which makes the study of the operators
$\BA_\a$ an interesting and non-trivial problem.

It turns out that the $\ba_0$-bound of the quadratic form $\bb[U]$ is
exactly $1/\sqrt 2$. This explains the role of the borderline value
$\a=\sqrt 2$. The techniques of quadratic forms does not apply to the large
values of $\a$. It was proved in \cite{nabsol} (and partly already in \cite{sol})
that the spectrum $\s(\BA_\a)$ for $\a>\sqrt 2$ is purely continuous and
coincides with
the whole of $\R$. The spectrum of the operator $\BA_{\sqrt 2}$ is also
purely continuous and coincides with the half-line $[0,\infty)$.

The operator $\BA_0$ can be easily studied by separation of variables.
It expands into the orthogonal sum of the operators in $\tL^2(\R)$ given
by
\begin{equation}\label{Hn}
\BH_n=-d^2/{dx^2}+(n+1/2),\qquad n\in\N_0:=\{0,1,\ldots\}.
\end{equation}
It follows that  the spectrum $\s(\BA_0)$ is absolutely continuous
and coincides with the half-line $[1/2,\infty)$. Its multiplicity
is $n$ on each interval $(n-1/2,n+1/2),\ n\in\N$.

The following statement, which describes the spectral properties
of  the operator $A_\alpha$ for  $\alpha$ small, is a particular
case of Theorem 6.2 in \cite{sol}.
\begin{prop}\label{Pr}
Let $\alpha<\sqrt 2$. Then
$\s_{\ess}(\BA_\a)=\s(\BA_0)=[1/2,\infty)$. The spectrum of
$\BA_\a$ below the threshold $\l_0=1/2$ lies in the interval
$(0,1/2)$, is always non-empty and consists of a finite number of
eigenvalues.
\end{prop}
This structure of the lower spectrum is typical for relatively
compact perturbations, a property which is violated in our case.
Here another mechanism is in effect and leads to the same result,
but only for small values of $\a$. We discuss this mechanism in
the final section 5.

\vskip0.3cm

Given a self-adjoint operator $\BT$ in a Hilbert space $\GH$
and a real number $s$, we denote
\begin{equation*}
N_+(s;\BT)=\dim E^\BT(s,\infty)\GH,\qquad N_-(s;\BT)=
\dim E^\BT(-\infty,s)\GH
\end{equation*}
where $E^\BT(\cdot)$ is the spectral measure of $\BT$. If, say,
$N_-(s;\BT)<\infty$, then the spectrum of the operator $\BT$ on
the interval $(-\infty,s)$ reduces to a finite number of
eigenvalues (counting their multiplicities), and $N_-(s;\BT)$ is
equal to this number.

\vskip0.3cm Our goal in this paper is study of the function
$N_-(1/2;\BA_\a)$ as $\a\nearrow\sqrt2$. We shall see that
$N_-(1/2;\BA_\a)\to\infty$, and calculate the asymptotics of this
function. This complements Proposition \ref{Pr} by giving a
quantitative characteristic of the discrete part of $\s(\BA_\a)$.
Probably, \thmref{t1} should be considered as the central result
of the paper. It establishes the equality
$N_-(1/2-\vare;\BA_\a)=N_+(\sqrt2/\a;\BJ(\vare))$ where
$\vare\in(0,1/2)$ and $\BJ(\vare)$ is a certain Jacobi matrix. For
the operator family $\BA_\a$ this is an analog of the classical
Birman -- Schwinger principle.

\section{Quadratic form  $\ba_\a$}

The quadratic form $\ba_0$ given by \eqref{ba0}
is positive definite and closed on the natural form-domain
\begin{equation*}
D:=\dom \ba_0=\{U\in\tH^1(\R^2): \ba_0[U]<\infty\}
\end{equation*}
where, as usual, $\tH^1$ stands for the Sobolev space.
The self-adjoint operator in $\tL^2(\R^2)$, generated by the
quadratic form $\ba_0[U]$, is $\BA_0$, i.e. the
operator \eqref{eq} -- \eqref{tran} for $\a=0$.

It is convenient to express both quadratic forms $\ba_0$ and $\bb$ in
terms of the decomposition of $U$ into the series in the (normalized in
$\tL^2(\R)$) Hermite functions in $y$:
\begin{equation}\label{U}
U(x,y)=\sum_{n\in\N_0}u_n(x)\chi_n(y).
\end{equation}
We often identify a function $U(x,y)$ with the sequence $\{u_n(x)\}$ and
write $U\sim\{u_n\}$. This identification is a unitary mapping of the space
$\tL^2(\R^2)$ onto
the Hilbert space $\ell^2(\N_0,\tL^2(\R))$. Let us recall the recurrence
relation for the functions $\chi_n$:
\begin{equation}\label{herm}
\sqrt{n+1}\chi_{n+1}(y)-\sqrt{2}y\chi_{n}(y)+\sqrt{n}\chi_{n-1}(y)=0,
\qquad n\in \N_0.
\end{equation}
Substituting in \eqref{ba0} the representation \eqref{U} of the function
$U$, we find
\begin{gather}
\ba_0[U]=\sum_{n\in\N_0}\bh_n[u_n]=:\sum_{n\in\N_0}\int_\R
\bigl(|u_n'|^2+(n+1/2)|u_n|^2\bigr)dx.\label{ba0s}
\end{gather}
The equality \eqref{ba0s} shows that the decomposition \eqref{U} diagonalyzes
the quadratic form $\ba_0[U]$ and hence, reduces the operator $\BA_0$. This
immediately implies the decomposition of $\BA_0$ into the orthogonal
sum of the operators $\BH_n$, see \eqref{Hn}, and hence
the structure of the spectrum $\s(\BA_0)$, described in the
Introduction.

In the same way, taking \eqref{herm} into account, we find that
\begin{equation}\label{bbs}
\bb[U]=\sum_{n\in\N_0}\sqrt{2n}\re\bigl(u_n(0)\overline{u_{n-1}(0)}\bigr).
\end{equation}

\bigskip
\begin{lem}\label{l1} (cf. \cite{sol}, section 6).
The quadratic form $\bb[U]$ is well-defined on $D$, and
\begin{equation}\label{fbound}
\sqrt2\,|\bb[U]|\le\ba_0[U],\qquad \forall U\in D.
\end{equation}
\end{lem}
\begin{proof}
Our argument is based upon the inequality
\begin{equation}\label{ineq}
2\g|u(0)|^2\le\int_\R\bigl(|u'|^2+\g^2|u|^2\bigr)dx,\qquad\forall
u\in\tH^1(\R),\ \g>0.
\end{equation}
Its proof is elementary and we skip it. It is also easy to show that the
equality in \eqref{ineq} is attained on the one-dimensional subspace in
$\tH^1(\R)$, generated by the function
\begin{equation}\label{at}
\wt u_\g(x):=(2\g)^{-1/2}e^{-\g|x|}.
\end{equation}
Here the factor $(2\g)^{-1/2}$ is chosen in such a way that
\begin{equation*}
\int_\R\bigl(|\wt u'_\g|^2+\g^2|\wt u_\g|^2\bigr)dx=1.
\end{equation*}

\vskip0.2cm

We derive from \eqref{bbs}:
\begin{gather*}
\sqrt2|\bb[U]|\le \sum_{n\in\N}\sqrt{n}(|u_n(0)|^2+|u_{n-1}(0)|^2)\\
=\sum_{n\in\N_0}\bigl(\sqrt{n}+\sqrt{n+1}\bigr)|u_n(0)|^2.
\end{gather*}
Since $\sqrt n +\sqrt{n+1}<\sqrt{2(2n+1)}$,
we conclude from \eqref{ineq} that
\begin{equation}\label{estbb}
\sqrt2|\bb[U]|\le\sum_{n\in\N_0}\bh_n[u_n]=\ba_0[U],\qquad \forall
U\sim\{u_n\}\in D,
\end{equation}
whence \eqref{fbound}.
\end{proof}

It follows from Lemma \ref{l1} that for $0<\a<\sqrt 2$ the quadratic
form $\ba_\a[U]$ is positive definite:
\begin{equation}\label{baa}
\ba_\a[U]\ge\bigl(1-\frac{\a}{\sqrt2}\bigr)\ba_0[U]\ge\frac1{2}
(1-\frac{\a}{\sqrt2}\bigr)\|U\|^2,\qquad U\in D
\end{equation}
(here and in the sequel $ \|U\|:=\|U\|_{\tL^2(\R^2)}$).
It is also closed, cf. e.g. \cite{birsol}, Lemma 1.1. The operator
$\BA_\a$ for such $\a$ can be defined as the self-adjoint operator
in $\tL^2(\R^2)$, associated with the quadratic
form $\ba_\a[U]$.

\section{Function $N_-(\frac1{2}-\vare;\BA_\a)$}

We are interested in the lower spectrum of the operator $\BA_\a$,
i.e. in the part of spectrum lying below the point
$1/2=\inf\s(\BA_0)$. By \eqref{baa}, this part of
$\s(\BA_\a)$ lies in the interval $1-\a/\sqrt2\le 2\l<1$. The
general perturbation theory gives no further information, since
the quadratic form $\bb$ is only relatively bounded but not
relatively compact with respect to $\ba_0$. However, we reduce the
problem to a simpler one, for a certain Jacobi operator in
$\ell^2(\N_0)$. This reduction allows us to handle the original
problem.

\vskip0.2cm
Fix a number $\vare,\ 0<\vare<1/2$ and consider a zero-diagonal Jacobi matrix
$\BJ(\vare)$ with the entries
\begin{equation*}
j_{n,n-1}(\vare)=j_{n-1,n}(\vare)=\frac{n^{1/2}}
{2(n+\vare)^{1/4}(n-1+\vare)^{1/4}},\qquad n\in\N.
\end{equation*}
All the other entries of the matrix are equal to zero. We use the same symbol
$\BJ(\vare)$ for the operator in $\ell^2(\N_0)$, generated by this matrix.
The operator $\BJ(\vare)$ is bounded and self-adjoint, its spectrum is invariant
under the reflection $\l\mapsto -\l$. It is well known that
$\s_{\ac}(\BJ(\vare))=[-1,1]$.
Besides, the operator may have (and actually, has) simple eigenvalues $\pm\l_n,
\l_n>1$, with the only possible accumulation points at $\l=\pm1$.

\begin{thm}\label{t1}
For any $\a\in(0,\sqrt 2)$, define $s(\a)=\sqrt 2/\a$. Then for
an arbitrary $\vare\in(0,1/2)$ the equality is satisfied:
\begin{equation}\label{equ}
N_-(1/2-\vare;\BA_\a)=N_+(s(\a),\BJ(\vare))=N_-(-s(\a),\BJ(\vare)).
\end{equation}
\end{thm}
The equality \eqref{equ} can be considered as one more
manifestation of the general Birman -- Schwinger principle.
\vskip0.2cm
\begin{proof}
According to the variational principle,
\begin{equation}\label{var}
N_-(1/2-\vare;\BA_\a)=\max_{\CF\in\GF(\vare)}\dim\CF
\end{equation}
where $\GF(\vare)$ is the set of all subspaces $\CF\subset D$, such that
\begin{equation}\label{var2}
\ba_\a[U]-(1/2-\vare)\|U\|^2_{\tL^2(\R^2)}<0,\qquad\forall U\in\CF,\
U\neq 0.
\end{equation}
Set
\begin{equation*}
\|U\|^2_\vare=\sum_{n\in\N_0}\int_\R\bigl(|u_n'|^2+(n+\vare)|u_n|^2\bigr)dx,
\qquad U\sim\{u_n\}.
\end{equation*}
For any $\vare>0$ and any $U\in\GH$ the quadratic form
$\|U\|^2_\vare$ can be estimated through $\ba_0[U]$ from above and
from below and hence, can be taken as a metric form on $D$. The
inequality \eqref{var2} can be re-written as
\begin{equation}\label{var3}
\|U\|^2_\vare+\a\sum_{n\in\N}\sqrt{2n}\re\bigl(u_n(0)\overline{u_{n-1}(0)}\bigr)<0.
\end{equation}

Consider the subspace $\wt{D}(\vare)$ in $D$, formed by the
elements
\begin{equation*}
\wt U\sim\{C_n\wt u_{\sqrt{n+\vare}}\},\ \{C_n\}\in\ell^2(\N_0),
\end{equation*}
where the elements $\wt u_\g$ for any $\g>0$ are given by \eqref{at}.
Note that
$\|\wt U\|_\vare=\|\{C_n\}\|_{\ell^2}$.
Let $\Pi_\vare$ stand for the operator which projects $D$ onto
$\wt{D}(\vare)$ and is orthogonal in the metric $\|\cdot\|_\vare$. If
$U\sim\{u_n\}\in D$, then
\begin{gather*}
\wt U_\vare:=\Pi_\vare U\sim\{C_n\wt u_{\sqrt{n+\vare}}\}
\end{gather*}
where
\begin{equation*}
C_n=\int_\R\bigl(u_n'\wt u'_{\sqrt{n+\vare}}+(n+\vare)u_n
\wt u_{\sqrt{n+\vare}}\bigr)dx=2^{1/2}(n+\vare)^{1/4}u_n(0).
\end{equation*}
We see that
\begin{equation}\label{rep}
C_n\wt u_{\sqrt{n+\vare}}(0)=u_n(0),\qquad \forall n\in\N_0.
\end{equation}
If in the inequality \eqref{var3} we replace $U$ by $\wt U_\vare$,
the first term in the left-hand side does not increase and the second remains unchanged,
so that the inequality remains valid. In other words, if a
subspace $\CF\subset D$ belongs to the class $\GF(\vare)$, then
also $\Pi_\vare\CF\in\GF(\vare)$.

On the other hand, assume that $\CF,\CF'$ are two subspaces of the
class $\GF(\vare)$, such that $\CF\subset\CF'$ and $\CF\subset \wt
D(\vare)$. Suppose that there exists an element $U\sim\{u_n\}\in\CF'$
orthogonal to $\CF$ in $\vare$-metric. Then by \eqref{rep} $u_n(0)=0$
for all $n\in\N_0$. This yields $\bb[U]=0$, which contradicts \eqref{var3}.
This means that our assumption implies $\CF=\CF'$.

It follows from these remarks that along with \eqref{var} the next
equality holds:
\begin{equation*}
N_-(1/2-\vare;\BA_\a)=\max_{\CF\in\GF(\vare),\ \CF\subset
\wt{D}(\vare)}\dim\CF.
\end{equation*}

For any $\wt U\sim\{C_n\wt
u_{\sqrt{n+\vare}}\}\in \wt{D}(\vare)$ we have
\begin{gather*}
\|\wt U\|^2_\vare+\a\bb[\wt U]
=\sum_{n\in\N_0}|C_n|^2+2s^{-1}\sum_{n\in\N}j_{n,n-1}(\vare)
\re(C_n\overline{C_{n-1}})\\
=\|g\|^2_{\ell^2}+s^{-1}\bigl(\BJ(\vare)g,g\bigr)_{\ell^2}, \qquad
g=\{C_n\}\in\ell^2.
\end{gather*}
The sum in the right-hand side is the quadratic form of the
operator $\BI+s^{-1}\BJ(\vare)$. Now \eqref{equ} is implied by the
variational principle and the symmetry of $\s(\BJ(\vare))$.
\end{proof}
\vskip0.2cm

\thmref{t1} does not apply to the most interesting case $\vare=0$,
since $j_{1,0}(0)=\infty$. However, we can restrict the quadratic
form $\bigl(\BJ(\vare)g,g\bigr)_{\ell^2}$ to the subspace
$\bigl\{g=\{C_n\}:\ C_0=0\bigr\}$ of codimension $1$. This may
shift the number of eigenvalues no more than by one. For the
problem obtained, the passage to the limit as $\vare\to 0$ is
already possible, and the resulting zero-diagonal Jacobi matrix is
$\BJ_0$ whose off-diagonal entries are given by
\begin{equation}\label{matr}
2j_{n,n-1}=2j_{n-1,n}=(1-n^{-1})^{-1/4},\qquad n-1\in\N.
\end{equation}

So, we arrive at the following result.
\begin{thm}\label{t2}
Let $\a\in(0,\sqrt 2)$ and $s(\a)=\sqrt 2/\a$. Then

\noindent
either $N_-(1/2;\BA_\a)=N_+(s;\BJ_0)$, or
$N_-(1/2;\BA_\a)=N_+(s;\BJ_0)+1$.
\end{thm}
\thmref{t2} reduces the problem of the asymptotic behaviour of the
function $N_-(1/2;\BA_\a)$ as $\a\nearrow\sqrt 2$ to the question
about the asymptotics of the eigenvalues of the
matrix $\BJ_0$, lying above the point $\l=1$.
We could not find the corresponding result in the literature,
so that we derive it the next section. Here is the formulation.

\begin{thm}\label{t3}
Let $\BJ$ be a zero-diagonal Jacobi matrix with the off-diagonal entries
\begin{equation}\label{ja}
j_{n,n-1}=j_{n-1,n}=1/2+qn^{-1}(1+o(1))
\end{equation}
where $q=const,\ q>0$. Then the operator $\BJ$ has the infinite number
of non-degenerate
eigenvalues $\pm \l_k(\BJ)$, such that
\begin{equation}\label{jaas}
\l_k(\BJ)=1+\frac{2q^2}{k^2}(1+o(1)),\qquad k\to\infty.
\end{equation}
These eigenvalues exhaust the spectrum of $\BJ$ outside $[-1,1]$.
Equivalently to \eqref{jaas},
\begin{equation}\label{jaas1}
N_+(s;\BJ)\sim\frac{q\sqrt2}{\sqrt{s-1}},\qquad s\searrow 1.
\end{equation}
\end{thm}

It follows from \eqref{matr} that the entries of the matrix $\BJ_0$
satisfy \eqref{ja} with $q=1/8$. Therefore, Theorems \ref{t2}
and \ref{t3} (equality \eqref{jaas1}) immediately imply the asymptotic formula
\begin{equation}\label{as}
N_-(1/2;\BA_\a)\sim\frac1{4\sqrt{2(s(\a)-1)}},\qquad
s(\a)=\sqrt2/\a,\ \a\searrow\sqrt 2.
\end{equation}

\section{Proof of \thmref{t3}}

The main ingredient of the proof is a result by W. Van Asshe
\cite{a} on a class of the orthogonal polynomials on the real
axis, namely of the so-called Pollaczek polynomials $P^\l(x;a,b)$,
see e.g. \cite{ch}. They depend on three real parameters $\l,a,b$
but we need only their particular case for $b=0,\ a=-r<0$ and
$\l>r$. The monic Pollaczek polynomials $Q^\l(x;r)$, i.e.
polynomials $P^\l(x;-r,0)$, normalized in such a way that their
leading coefficient becomes $1$, satisfy the recurrent relation
\begin{gather*}
Q^\l_{n+1}(x;r)=xQ^\l_{n}(x;r)-p_n(\l,r)
Q^\l_{n-1}(x;r),\\
p_n(\l,r)=\frac{n(n+2\l-1)}{4(n-r+\l-1)(n-r+\l)},\qquad n\in\N.
\end{gather*}
The polynomials $Q^\l_{n}$ correspond to the zero-diagonal Jacobi matrix $\BJ(\l,r)$
whose off-diagonal entries are
\begin{equation}\label{entr}
j_{n,n-1}=j_{n-1,n}=\sqrt{p_n(\l,r)}.
\end{equation}
It was proven in \cite{a}, Section III that the spectrum of $\BJ(\l,r)$,
lying outside the segment $[-1,1]$,
consists of the non-degenerate eigenvalues $\pm\mu_k=\pm\mu_k(\l,r)$ where
$\mu_k$ satisfy the equation
\begin{gather*}
\l-\frac{r\mu}{\sqrt{\mu^2-1}}=-k,\qquad k\in\N_0.
\end{gather*}
This gives
\begin{equation}\label{mu}
\mu_k=\biggl(1-\frac{r^2}{(k+\l)^2}\biggr)^{-1/2}=1+\frac{r^2}{2k^2}
+o\bigl(\frac1{k^2}\bigr),\qquad k\to\infty.
\end{equation}

Another ingredient is a variational principle for the eigenvalues of
Jacobi matrices, see \cite{g}, Lemma III.1. This variational principle is
an almost immediate consequence of Sturm's comparison theorem, see e.g.
\cite{f}, Theorem 1 on p. 152. Below we present its formulation for a
particular case we need in this paper.

\begin{lem}\label{le}
Let $\BJ,\ \BJ'$ be Jacobi matrices with the zero diagonal entries
and the off-diagonal entries $j_{n,n-1}=1/2+b_n, \
j'_{n,n-1}=1/2+b'_n$, such that $0\le b_n\le b'_n$ for all $n$
and $b_n'\to 0$. Then $\s_{ess}(\BJ)=\s_{ess}(\BJ')=[-1,1]$ and
for any $s>1$
\begin{equation*}
N^+(s;\BJ)\le N^+(s;\BJ').
\end{equation*}
\end{lem}

\bigskip
Now we are in a position to prove \thmref{t3}. It follows from
Lemma \ref{le} that the asymptotic behaviour of the eigenvalues
does not depend on the term $o(1)$ in \eqref{ja}. The entries
$j_{n,n-1}$ in \eqref{entr} satisfy $2j_{n,n-1}\sim 1+r/n$.
Clearly, \eqref{jaas} is a direct consequence of \eqref{mu}.
\vskip0.2cm

\section{Concluding remarks}
{\bf 5.1.} Here we explain, why for $\a<\sqrt2$ the structure of
the spectrum $\s(\BA_\a)$ of the operator $\BA_\a$ below the
threshold $1/2$ is the same as if the perturbation were relatively
compact. Of course, one such explanation is given by the proof of
\thmref{t1}, however we shall present also one more argument, of a
somewhat more heuristic nature. A rigorous version of this
argument was used in \cite{sol} for the proof of Theorem 6.2.

The quadratic form $\bb[U]$, see \eqref{bb} and \eqref{bbs}, is
the sum of terms
\begin{equation*}
\bb_n[U]=\sqrt{2n}\re\bigl(u_n(0)\overline{u_{n-1}(0)}\bigr),
\end{equation*}
each of rank two. The quadratic form $\bb_n[U]$ interacts only
with the terms $\bh_{n-1}[u_{n-1}]$ and $\bh_n[u_n]$ in the
representation \eqref{ba0s} of the quadratic form $\ba_0[U]$. The
term $\bh_n[u_n]$ corresponds to the operator $\BH_n$, see
\eqref{Hn}, whose spectrum is $[n+1/2,\infty)$. The perturbation
of the spectrum, brought by the term $\a\bb_n[U]$, does not reach
the point $\l_0=1/2$, provided that $\a<\sqrt 2$ and $n$ is large
enough. This means that effectively we are dealing with a finite
rank perturbation, as soon as we restrict ourselves with the small
values of the coupling parameter and are interested only in the
lower part of $\s(\BA_\a)$. \vskip0.2cm
{\bf 5.2.} In the
paper \cite{sol} the operator family $\BA_\a$ was considered in a
more general setting. Namely the operators act in the space
$\tL^2(\G\times\R)$ where $\G$ is a {\sl metric star graph}, i.e.
a graph with $m$ bonds $\CB_1,\ldots,\CB_m, \ 1\le m<\infty$, all
emanating from a common vertex $o$. Let us recall that each bond
of a metric graph is viewed as a line segment of finite or
infinite length. The real axis $\R$ can be considered as the star
graph with two bonds (so that $m=2$), each of infinite length, and
with $o=0$.

Let us identify each bond $\CB_j$ with the segment $[0,B_j)$,
where $B_j\le\infty$ is the length of $\CB_j$. We denote by
$x$ the coordinate along each bond (dropping the index $j$); the
value $x=0$ corresponds to the point $o$.

The action of the operator $\BA_\a$ in this, more general case is
defined by the same equality \eqref{eq}, in which $y$ denotes the
coordinate along the additional straight line. The
condition \eqref{tran} is replaced by the matching conditions
\begin{gather*}
U^{1}(0,y)=\ldots=U^{m}(0,y);\\
U^{1}_x(0,y)+\ldots+U^{m}_x(0,y)=\a y U(0,y)
\end{gather*}
where $U^{j}$ stands for the restriction of $U$ to the bond $\CB_j$. Besides,
the Dirichlet condition $U^{j}(B_j,y)=0$ is imposed for each bond of finite length.

Theorem 6.2 in \cite{sol} (cf. Proposition \ref{Pr} of the present
paper) was proved for this general version of the operator
$\BA_\a$. The only difference with the particular case $\G=\R$ is
that the $\ba_0$-bound of the quadratic form $\bb$ is $\sqrt 2/m$,
cf. \eqref{fbound}. Correspondingly, the techniques developed in
the present paper allows one to prove an analog of the asymptotic
relation \eqref{as}. The only distinction is that for any star
graph with $m$ bonds we have to take $\a\nearrow m/\sqrt 2$ and
$s(\a)=m/(\a\sqrt2)$.

\bigskip

In conclusion, I would like to express my deep gratitude to S.N.
Naboko for very useful discussions, and to the referee for
pointing out an arithmetic error in calculation of the asymptotic
coefficients. This error is corrected in the final version of the
paper.

\bibliographystyle{amsplain}

\end{document}